\documentclass[smallcondensed,final]{svjour3}

\journalname{Acta. Applicandae Mathematicae}

\usepackage{epsfig}
\usepackage{psfig}
\usepackage{graphicx}
\usepackage{amsmath}
\usepackage{amssymb}
\usepackage{geometry}
\everymath{\displaystyle}

\usepackage{xcolor}

\newtheorem{rem}{\bf Remark}

\newcommand{\ede}{\hfill \rule [-.3em]{.5em}{.5em}}

\newcommand{\reel}{I\!\!R}

\newcommand{\beqnum}{\begin{equation}\begin{array}{lcl}}
\newcommand{\eeqnum}{\end{array}\end{equation}}
\newcommand{\beqnom}{\begin{eqnarray}}
\newcommand{\eeqnom}{\end{eqnarray}}
\newcommand{\beqnc}{\begin{center}\begin{eqnarray}}
\newcommand{\eeqnc}{\end{eqnarray}\end{center}}
\newcommand{\beqnlm}{\begin{equation}\vspace{-.5cm}\begin{array}{lll}}
\newcommand{\eeqnlm}{\end{array}\end{equation}}\vspace{-.5cm}
\newcommand{\beq}{\begin{eqnarray*}}
\newcommand{\eeq}{\end{eqnarray*}}
\newcommand{\bef}{\begin{figure}}
\newcommand{\enf}{\end{figure}}

\usepackage{epsfig}
\usepackage{mathptm}

\usepackage{times}
\begin{document}

\title{Path following for a target point attached to a unicycle type vehicle}
\author{S. Laghrouche  \and Y. Chitour \and
M. Harmouche \and
F. S. Ahmed}
\institute{S. Laghrouche \and M. Harmouche \and F. S. Ahmed  \at Laboratoire Systemes et Transports, Universite de Technologie Belfort-Montbeliard, Belfort, France \email{[salah.laghroughe;mohamed.harmouche;fayez-shakil.ahmed]@utbm.fr}
\and
Y. Chitour \at L2S, Université Paris-Sud XI, CNRS and Supélec, Gif- sur-Yvette, France \email{yacine.chitour@lss.supelec.fr}
}

\maketitle

\begin{abstract}                          
In this article, we address the control problem of unicycle path following, using a rigidly attached target point. The initial path following problem has been transformed into a reference trajectory following problem, using saturated control laws and a geometric characterization hypothesis, which links the curvature of the path to be followed with the target point. The proposed controller allows global stabilization without restrictions on initial conditions. The effectiveness of this controller is illustrated through simulations.
\end{abstract}

\section{Introduction}

The case of vehicle path following using a "target point" (situated at a distance from the vehicle) is well known in the domain of automatic vehicle guidance. This technique is often used in robotic vehicles with artificial camera vision, where the camera is fixed on the vehicle and the target point (physical or virtual) is situated somewhere in its field of view. This problem has been the subject of many research works in the recent years \cite{consolini1,consolini2,Jiang,Kanayama,MS,Samson,SamsonAit}. The dominant trend in the contemporary literature is to control either the vehicle's forward velocity (thereby, not controlling the vehicle's orientation), or the instantaneous rotational velocity only. Hence, essentially only one actuator is used.

In \cite{consolini2}, a local path following strategy has been proposed, which takes uncertainties into account as well. Their solution is based on a control law that comprises of two terms; an open loop control that allows inversion of the nominal model, and a closed loop control that stabilizes the resulting system. It should be noted that the error dynamics obtained in \cite{consolini2} are expressed in the Fr\'{e}net frame associated to the followed path (a technique that has also been discussed in \cite{MS}). While the use of Fr\'{e}net frames is convenient, its application is \emph{local}, i.e. the convenience is significant only when the vehicle is close to the path (with respect to a universal constant), positioned and oriented. When such ideal situations are not present and the vehicle is actually located far from the path, another controller (e.g. an open loop control) takes over to bring the vehicle in the path's proximity before the primary controller starts operation.

In \cite{Makarov1,Makarov2}, a polar state transformation has been used to ensure tracking of smooth plane trajectory for a trailer-truck vehicle. This coordinate transformation is not global, therefore the proposed controllers only ensure local stabilization of the system.
In \cite{Egerstedt}, controllers have been proposed to follow a reference trajectory using a virtual vehicle approach, where the motion of the reference point on the desired trajectory is governed by a differential equation containing the error feedback. The proposed controllers do not ensure convergence to the reference trajectory itself. In fact the convergence is limited to a bounded distance from the trajectory. This is due to the fact that this controller stabilizes the angle variable exponentially.

In this paper, we present a target point based path following technique for a robot unicycle. The target point has been considered fixed with respect to a point on the vehicle. More precisely, the target point is at a fixed distance $d>0$ from the center of gravity on the axis of the vehicle. Our control objective is to drive the vehicle, such that the target point follows the desired path (see Fig. 1 below). We have assumed that the vehicle's velocity is measured only, and not controlled. This assumption conforms with practical applications,
where other intelligent systems control the velocity, (for example, ABS, ESP \cite{pasillias}).

The primary objective of this work is to conceive \emph{global} control laws, which are applicable regardless of the initial position and orientation of the vehicle w.r.t. the path to be followed. Hence the problem can be defined as orientation control with a forced forward velocity. Our solution is based upon parameterization of the reference path as ''the trajectory of a unicycle'', the forward velocity of which can be considered as a supplementary control variable. A similar approach can be found in  \cite{consolini1} and \cite{consolini2}, where orientation control of a vehicle is under consideration. The authors have achieved this through a dynamic inversion process, implemented using adaptive parametrization of the followed path. In our work, we have chosen the opposite direction, converting the problem of path following into a special case of trajectory following. Furthermore, we have also considered the trajectory of the target point as the trajectory of a unicycle. This allows us to express the error dynamics as the difference between the unicycle dynamics defined by the reference path, and the unicycle dynamics defined by the target point. We have thus obtained a controlled system with three dimensional state and two control inputs (the forward velocity of the reference path and the angular velocity of the vehicle) .

Our control law is based upon state feedback with static error control algorithms, along with saturated input technique
\cite{Yakoubi2006,Yakoubi2007,Yakoubi2007-1,Yakoubi2007-2,teel,chitour1,chitour2,chitour3}. As would be shown further on, the application of bounded inputs is justified by two constraints, $(a)$ to maintain the forward velocity on the reference path uniformly bounded, ($b)$ to focus on controlling the orientation of the unicycle defined by the target point, rather than controlling the orientation of the vehicle. It is worth mentioning that in order to satisfy constraint $(b)$, we have supposed the geodesic curvature of the followed path to be strictly bounded in magnitude by the inverse of the distance $d$. Application of such type of bounded commands in the same context  (trajectory following of unicycle robots) can be found in \cite{Jiang} . The stability analysis is based on an argument of the Lyapunov type. Our contribution, compared to \cite{Jiang} is the determination of a \emph{strict and global} Lyapunov Function on an appropriate basin of attraction. As a byproduct, we can handle  model uncertainties, external perturbations as well as (constant) delays as indicated in a series of remarks preceding the simulation section.


\vspace{2mm}

{\bf Acknowledgements.} The authors thank E. Panteley and W. Pasillas-L\'epine for their constructive remarks.
\section{Vehicle model and reference trajectory}
Let us consider a path $\gamma$ with geodesic curvature  $\kappa_r^*$ whose absolute value is bounded by $\kappa_{max}\geq 0$. As described in the introduction, we want to parameterize $\gamma$ as a unicycle trajectory with a forward velocity $u(t)$ such that $\gamma(t)=(p_r(t),q_r(t))$ can be described by the following state equations :
\begin{eqnarray}
\begin{array}{rcl}
\dot{p_r}&=&u\cos{\psi_r},\\
\dot{q_r}&=&u\sin{\psi_r},\\
\dot{\psi_r}&=&u\kappa_r,
\end{array}
\end{eqnarray}

where $\kappa_r$, is the scalar curvature associated to the parametrization of $\gamma$ by time $t$. The relationship between the arclength $s$ of $\gamma$ and time $t$ for the trajectory $(p_r,q_r,\psi_r)$ is given by $s(t)=s_0+\int_0^tu(\tau)d\tau$. The scalar curvature $\kappa_r(t)$ is hence equal to  $\kappa_r^*(s(t))$. For the sake of simplicity, we have assumed in this paper that $u$ is a strictly positive function (i.e., strictly positive forward velocity), and moreover, that the controls $u$ verify $\int_{0}^{\infty}u(t)dt=+\infty$. Furthermore, for all $t\geq 0$, we have
\begin{equation}
|\kappa_r(t)|\leq \kappa_{max}.
\end{equation}

The state equations for the vehicle can be defined as:
\begin{eqnarray}
\begin{array}{rcl}
\dot{x}&=&V_x\ \cos{\psi},\\
\dot{y}&=&V_x\ \sin{\psi},\\
\dot{\psi}&=&V_x\  v,
\end{array}
\end{eqnarray}
These equations represent the vehicle's motion with a velocity $V_x$, along the curve defined by its geodesic curvature $v$. This variable will be considered as the second control in the problem. Notice that $V_x$ is not necessarily constant, but simply a continuous function of time, which verifies the following hypothesis: there exist two positive constants $0<V_{min}\leq V_{max}$, such that for all $t\geq 0$
\begin{equation}\label{vx0}
V_{min}\leq V_x(t)\leq V_{max}.
\end{equation}
Recall that the strict positivity of the lower bound is a necessary assumption to obtain the results of the paper (see, \cite{Samson} for an explanation of this classical phenomenon). Indeed, Eq. \ref{vx0} implies that linearized systems associated to the reference trajectory are controllable and thus the nonlinear system is locally controllable.

For the target point, the equations for the coordinates $p$ and $q$ are defined as:
\begin{eqnarray}
\begin{array}{rcl}
p&=x+ d \cos{\psi},\\
q&=y+ d \sin{\psi}.\\
\end{array}
\end{eqnarray}

We will also suppose throughout the paper that
\begin{itemize}
\item [{\textbf (H1)}] $d\kappa_{max}<1$.
\end{itemize}

This can be considered as a technical condition, or a design constraint for positioning the target point. However, as explained later, condition $(H1)$ turns out to be (almost) necessary to control the system.

The dynamics of the target point can be obtained by deriving the precedent equations
\begin{eqnarray}
\begin{array}{rcl}
\dot{p}&=&V_x\ \cos{\psi}  - d\ V_x\ \sin{\psi}\ v, \\
\dot{q}&=&V_x\ \sin{\psi} + d\ V_x\ \cos{\psi}\ v, \\
\dot{\psi}&=&V_x\  v.
\end{array}
\end{eqnarray}

\begin{figure}[ht!]
\begin{center}
  \includegraphics[width=8cm]{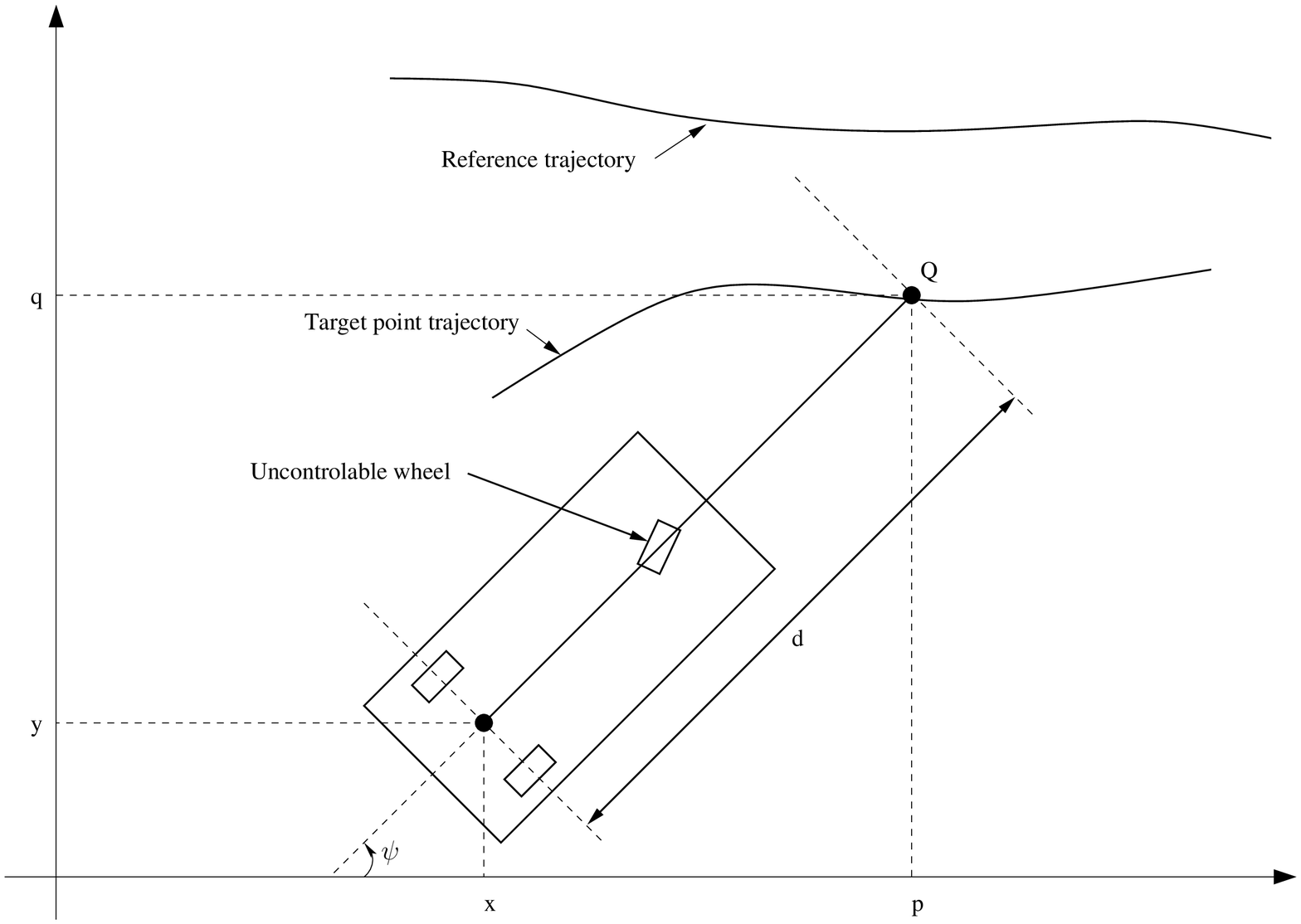}
  \caption{The reference trajectory, the vehicle and its target point.}
\end{center}

\end{figure}

The curve defined by the target point is traveled at the following speed:
$$
v_d:= \sqrt{\dot p^2+\dot q^2}=V_x \sqrt{1+(vd)^2}.
$$

Our objective now is to define the dynamics of the target point as those of a unicycle. Therefore, let us consider $\theta$ as the angle between the abscissa axis and the velocity vector $(\dot p,\dot q)^T$. It can easily be seen that $\theta =\psi +\arctan(d \nu)$,  and therefore,
$$
\dot p= v_d\cos(\theta),\ \ \ \dot q=v_d\sin(\theta),
$$

The scalar curvature $\omega$ is defined by
${\displaystyle\omega:=\frac{\dot{\theta}}{v_d}}$.

Solving these equations, we obtain:

\begin{equation}
{\displaystyle\omega=\frac{V_x\ v}{v_d}+\frac{d\ \dot{v}}{v_d (1+(v\
d)^2)}}.\label{omega}
\end{equation}

Hence the dynamics of the target point $(p,q)$ becomes
\begin{eqnarray}
\begin{array}{rcl}
\dot{p}&=&v_d\ \cos{\theta},\\
\dot{q}&=&v_d\ \sin{\theta},\\
\dot{\theta}&=&v_d\  \omega.
\end{array}
\end{eqnarray}

From here on, we will replace $v$ with $\omega$ as the new control. Considering equation (\ref{omega}), we obtain the following form:

\begin{equation}\label{edo-v}
{\displaystyle\dot{v}=\frac{1 + (v d)^2}{d} V_x \left[\sqrt{1 + (v
d)^2} \omega - v\right]},
\end{equation}

i.e. an ordinary differential equation for the unknown function $v$. Since the right side of (\ref{edo-v}) is not globally Lipschitz with respect to $v$, the solution may only be defined for finite time duration. We will show later on, that a choice of $\omega$ under  Hypothesis {\textbf (H1)} solves this problem (see Lemma~\ref{explosion} below).

The error between the target point and the reference curve can be defined as:
\begin{eqnarray}
\begin{array}{rcl}
e_p&=&p-p_r,\\
e_q&=&q-q_r,\\
\xi&=&\theta - \psi_r,
\end{array}
\end{eqnarray}

and the error dynamics are given by:
\begin{eqnarray}\label{erreur}
\begin{array}{rcl}
\dot{e_p}&=&v_d\ \cos{\theta} - u\cos{\psi_r},\\
\dot{e_q}&=&v_d\ \sin{\theta} - u\sin{\psi_r},\\
\dot{\xi}&=&v_d\  \omega - \kappa_r\ u.
\end{array}
\end{eqnarray}

The objective, hence, is to determine the control laws, $u(t,e_p,e_q,\xi)$
and $\omega(t,e_p,e_q,\xi)$ such that the closed loop system (\ref{erreur}) is globally asymptotically stable (GAS for short) with respect to the origin.

Let us first of all perform a variable change on the control, as follows:
\begin{eqnarray}\label{com0}
\begin{array}{rcl}
u&=&v_d(1+u_1),\\
\omega&=&\kappa_r(1+u_1) + u_2.
\end{array}
\end{eqnarray}
The system that we have to stabilize, becomes:
\begin{eqnarray}\label{err1}
\begin{array}{rcl}
\dot{e_p}&=&v_d(\cos{\theta}- \cos{\psi_r}-\ u_1\cos{\psi_r}),\\
\dot{e_q}&=&v_d(\sin{\theta} - \sin{\psi_r}-\ u_1\sin{\psi_r}),\\
\dot{\xi}&=&v_du_2.
\end{array}
\end{eqnarray}

The following lemma  provides bounding conditions on  $u_1$ and $u_2$ that would guarantee that the differential equation given in (\ref{edo-v}) is
defined for all times $t\geq 0$.

\begin{lemma}\label{explosion}
Suppose that for all $t\geq 0$, there exists
\begin{equation}\label{cond0}
{\displaystyle\frac{|u_1(t)|}d+|u_2(t)|\leq \beta_M
:=\frac{1-d\kappa_{max}}d}.
\end{equation}
Then, the differential equation given by Eq.(\ref{edo-v}) is defined for all times $t\geq 0$.
\end{lemma}
{\bf Proof of Lemma \ref{explosion}}
 Let us multiply
(\ref{edo-v}) by $v$. We obtain:
\begin{equation}
{\displaystyle v\dot{v}=\frac{1 + (v d)^2}{d} V_x \left[\sqrt{1 + (v
d)^2} v\omega - v^2\right]}.
\end{equation}

If $v\omega\leq 0$, then $v\dot v\leq 0$. If $v\omega>0$, the precedent equation can be written as:
\begin{equation}\label{equadiff2}
{\displaystyle  v\dot{v}=\frac{1 + (v d)^2}{d} V_x\ |v|
\left[\frac{\omega^2 + v^2 \left((d\ \omega)^2-1\right)}{\sqrt{1 +
(v d)^2} |\omega| + |v|}\right]}.
\end{equation}
In order to guarantee that the right side of (\ref{equadiff2}) is globally Lipschitz with respect to $v^2$, it is sufficient to choose $u_1,u_2$ such that for all $t\geq 0$,
$$
(d\ \omega(t))^2-1\leq 0.
$$
Using (\ref{com0}), we can rewrite the equation
$$
(d\ \omega(t))^2 -1\leq d^2(\kappa_{max}(1+|u_1(t)|)+|u_2(t)|)^2-1.
$$
For the value of this quantity to be less than zero, it is sufficient that $d(\kappa_{max}(1+|u_1(t)|)+|u_2(t)|)\leq 1$, and hence for all $t\geq 0$,
$$
{\displaystyle \frac{|u_1(t)|}d+|u_2(t)|\leq \beta_M.}
$$

\ede

In order to verify (\ref{cond0}), the controls $u_1$ et $u_2$ can be expressed in the following form:
\begin{eqnarray}
\begin{array}{rcl}
u_1&=&C_1\sigma(\cdot),\\
u_2&=&\beta \sigma(\cdot),
\end{array}
\end{eqnarray}
with (for instance)
\begin{equation}\label{eee}
(Cond0)\quad\quad {\displaystyle 0<C_1\leq \frac{d\beta_M}{2}}, \ \
{\displaystyle 0<\beta\leq \frac{\beta_M}2},
\end{equation}
and $\sigma$ being equal to the standard saturation function
$$
{\displaystyle \sigma(x)=\frac{x}{max(1,|x|)}}.
$$

Since $v$ is bounded, $v_d$ also remains uniformly bounded throughout $t\geq 0$. We can hence change the time by considering $dt'=v_d\ dt$. To keep the notations simple, we would continue to use $t$ for time. This has no effect on the control laws since our design is based on static feedback (w.r.t. the error).

The error dynamics hence becomes:
\begin{eqnarray}
\begin{array}{rcl}
\dot{e_p}&=&\cos{\theta}- \cos{\psi_r}-u_1\cos{\psi_r}, \\
\dot{e_q}&=&\sin{\theta} - \sin{\psi_r}-u_1\sin{\psi_r}, \\
\dot{\xi}&=&u_2.
\end{array}
\end{eqnarray}

Let us perform the following change of variable corresponding to a time-varying
rotation in the frame of the reference trajectory:

\begin{eqnarray}
\begin{array}{rcl}
y_1&=&e_p\ \cos{\psi_r} + e_q\ \sin{\psi_r}, \\
y_2&=&-e_p\ \sin{\psi_r} + e_q\ \cos{\psi_r}.
\end{array}
\end{eqnarray}

The final system can be expressed as

\begin{eqnarray}\label{system}
\begin{array}{rcl}
\dot{y_1}&=&-u_1 + (\cos{\xi}-1) + (1+u_1) \kappa_r y_2,\\
\dot{y_2}&=&\sin{\xi} - (1+u_1) \kappa_r y_1,\\
\dot{\xi}&=&u_2.
\end{array}
\end{eqnarray}

This system of equations greatly resembles the error dynamics obtained for the classic tracking problem of a vehicle using a unicycle, with the forward velocity and the instantaneous rotation velocity of the vehicle body as control variables (cf. \cite{MS} et \cite{Jiang}).

We choose the controls $u_1$ and $u_2$ as follows:
\begin{eqnarray}\label{oup0}
\begin{array}{rcl}
u_1&=&C_1 \sigma(M y_1) ,\\
u_2&=&\beta \sigma({\displaystyle\frac {- C_0}{\beta}\left[\xi + \rho
\sigma(C_2 y_2)\right]}),
\end{array}
\end{eqnarray}

with $M, C_0,C_1, C_2,\beta, \rho$ as positive constants to be fixed later.

Hence the error dynamics are:
\begin{eqnarray}\label{system1}
\begin{array}{rcl}
\dot{y_1}&=&-C_1 \sigma(M y_1)+ \lambda(t) y_2 + (\cos \xi-1),\\
\dot{y_2}&=&\sin{\xi} - \lambda(t) y_1,\\
\dot{\xi}&=&-\beta \sigma(\frac{C_0}{\beta}\left[\xi + \rho
\sigma(C_2 y_2)\right]),
\end{array}
\end{eqnarray}

where $\lambda(t) := (1+u_1) \kappa_r$. $\lambda$ is bounded by

\begin{eqnarray}\label{lam0}
\begin{array}{rcl}
|\lambda(t)|\leq (3+C_1) \kappa_{max}.
\end{array}
\end{eqnarray}

\textbf{Theorem 1:}
With the control $u_1$ and $u_2$ defined in \eqref{oup0}, for all $C_1$ and $\beta$ verifying (\ref{eee}), the system (\ref{system1}) is GAS with respect to $0$.

{\bf Proof of Theorem 1}

We first have the following result, which is a trivial consequence of the dynamics of $\xi(\cdot)$.
\begin{lemma}\label{lemt_0}
For every trajectory of \eqref{system1}, there exists $t_0 \geq 0$ such that, for every $t > t_0$ : $|\xi(t)| < 2 \rho$.
\end{lemma}
This follows from the fact that if $\left| \xi  \right| > \frac{3}{2}\rho $, then $\left[ {\xi  + \rho \sigma \left( {{C_2}{y_2}} \right)} \right]$ has the same sign of $\xi$, and $\left| {\xi  + \rho \sigma \left( {{C_2}{y_2}} \right)} \right| > \frac{\rho }{2}$. Finally, we get $\xi \dot \xi  <  -\frac{\rho}{2} \beta \sigma \left( {\frac{{{C_0}\rho }}{{2\beta }}} \right) < 0.$

We next impose  the following condition.
$$
(Cond1): \quad \quad 3\rho C_0\leq \beta.
$$

This implies that for $t\geq t_0$,

$$
\left|\frac{C_0}{\beta}\left[\xi(t)+ \rho \sigma(C_2
y_2(t))\right]\right|\leq 1.
$$

Hence, for $t\geq t_0$, the system (\ref{system1}) becomes:
\begin{eqnarray}\label{system2}
\begin{array}{rcl}
\dot{y_1}&=&-C_1 \sigma(M y_1)+ \lambda(t) y_2 + (\cos \xi -1),\\
\dot{y_2}&=&\sin{\xi} - \lambda(t) y_1,\\
\dot{\xi}&=&-C_0\left[\xi + \rho \sigma(C_2 y_2)\right].
\end{array}
\end{eqnarray}

Let $\cal{E}$ be a set of points $(y_1,y_2,\xi)$ such that $|\xi|<2 \rho$. According to Lemma \ref{lemt_0}, $\cal{E}$ is an open invariant set for the system (\ref{system2}). To prove Theorem 1, it is sufficient to form a strict Lyapunov function for (\ref{system2}) on $\cal{E}$. We propose the following candidate function:

\begin{eqnarray}
\begin{array}{rcl}
V(y_1,y_2,\xi):= {\displaystyle \frac{y_1^2+y_2^2}{2}+ \frac{F(\xi)
y_2}{C_0}+ \frac{N}{2C_0}\xi^2},
\end{array}
\end{eqnarray}
with $N$ a positive constant to be determined, and
$F(\xi)=\displaystyle{\int^{\xi}_0{{\displaystyle\frac{\sin{s}\ ds}{s}}}}$.


Notice that $F$ is an odd function, and if $N>{\displaystyle\frac1{C_0}}$ then $V$ is positive definite. We next prove that $V$ is a strict Lyapunov function for (\ref{system2}) on $\cal{E}$ with an appropriate choice of the constants.


Let us suppose from this point on that $\rho  \leqslant \frac{1}{2}$. Therefore, for $\vert \xi \vert\leqslant 2\rho $, one has
\beqnum\label{demo_1}
\begin{gathered}
  1 - {\frac{\xi }{6}^2} \leqslant \frac{{\sin \xi }}{\xi } \leqslant 1 ,\hfill \\
  1 - {\frac{\xi }{{18}}^2} \leqslant \frac{{F(\xi )}}{\xi } \leqslant 1 ,\hfill \\
	1 - \frac{{{\xi ^2}}}{2} \leqslant \cos \xi  \leqslant 1 .\hfill \\
\end{gathered}
\eeqnum
From here, it can be deduced that:
\beqnum\label{demo_2}
\begin{gathered}
  1 - {\frac{{2\rho }}{3}^2} \leqslant \frac{{\sin \xi }}{\xi } \leqslant 1, \hfill \\
  1 - {\frac{{2\rho }}{9}^2} \leqslant \frac{{F(\xi )}}{\xi } \leqslant 1, \hfill \\
  1 - 2{\rho ^2} \leqslant \cos \xi  \leqslant 1 .\hfill\\
\end{gathered}
\eeqnum
The derivative of $V$ along the trajectories of the system is equal to:

\beqnum\label{demo_3}
\dot V &=&  - \left[ {{C_1}{y_1}\sigma (M{y_1}) - \frac{{\lambda (t)F(\xi )}}{{{C_O}\xi }}\xi {y_1} + \frac{1}{2}\left( {N - \frac{{F(\xi )\sin \xi }}{{{C_O}{\xi ^2}}}} \right){\xi ^2}} +y_1( cos \xi -1 )  \right]\\
&&- \left[ {\frac{1}{2}\left( {N - \frac{{F(\xi )\sin \xi }}{{{C_O}{\xi ^2}}}} \right){\xi ^2} + \rho N\xi \sigma \left( {{C_2}{y_2}} \right)+ \frac{{\sin \xi }}{\xi }\rho {y_2}\sigma \left( {{C_2}{y_2}} \right)} \right].
\eeqnum

From equations (\ref{lam0}) and (\ref{demo_2}), it can be seen that the first term in brackets of equation (\ref{demo_3}) is greater or equal to:
\beqnum\label{demo_4}
A\left( {y{}_1,\xi } \right)&:=& {C_1}{y_1}\sigma \left( {M{y_1}} \right) - \frac{{\left( {3 + {C_1}} \right){\kappa _{max}}}}{{{C_O}}}\left| {\xi {y_1}} \right|\\
&& - \frac{1}{2}{\xi ^2}\left| {{y_1}} \right| + \frac{1}{2}\left( {N - \frac{1}{{{C_O}}}} \right){\xi ^2}.
\eeqnum
Similarly, the second term  in brackets of equation (\ref{demo_3}) can be bounded by:
\beqnum\label{demo_5}
B\left( {y{}_2,\xi } \right)&:=& \frac{1}{2}\left( {N - \frac{1}{{{C_O}}}} \right){\xi ^2} - \rho N\left| {\xi \sigma \left( {{C_2}{y_2}} \right)} \right| \\
&&+ \left( {1 - \frac{{2{\rho ^2}}}{3}} \right)\rho {y_2}\sigma \left( {{C_2}{y_2}} \right).
\eeqnum
Hence, using equations (\ref{demo_3}), (\ref{demo_4}) and (\ref{demo_5}), $\dot{V}$ can be expressed as:
\beqnum\label{demo_6}
\dot V \leqslant  - A\left( {y{}_1,\xi } \right) - B\left( {y{}_2,\xi } \right).
\eeqnum
We shall now present two lemmas, and establish the conditions on constants, under which these lemmas would hold true.
\begin{lemma}\label{lem_demo1}
\emph{There exist constants $C_1,\rho,\beta,M,N,C_0$ for which the function} $A$ \emph{is positive definite on} $\reel \times \left] { - 2\rho ,2\rho } \right[$.\\
\end{lemma}

\begin{lemma}\label{lem_demo2}
\emph{There exist constants $C_1,\rho,\beta,M,N,C_0$ for which the function }$B$ \emph{is positive definite on} $\reel \times \left] { - 2\rho ,2\rho } \right[$.\\
\end{lemma}

\emph{Proof of Lemma \ref{lem_demo1}:} Let us consider 2 cases:\\
\textbf{Case 1:}  $\left| {{y_1}} \right| \geqslant \frac{1}{M}$.\\
As $\left| \xi  \right| \leqslant 2\rho $, we obtain:
\beqnum\label{demo_7}
A \geqslant \left| {{y_1}} \right|\left( {{C_1} - \frac{{2\rho {\kappa _{max}}}}{{{C_0}}}\left( {3 + {C_1}} \right)  -{2{\rho ^2}} } \right).
\eeqnum
Hence it is sufficient to verify that:
\beqnum\label{demo_8}
{C_1} - \frac{{2\rho {\kappa _{max}}}}{{{C_0}}}\left( {3 + {C_1}} \right) -{2{\rho ^2}} > 0 \Leftrightarrow {C_1}\left( {1 - \frac{{2\rho {\kappa _{max}}}}{{{C_0}}}} \right) > \frac{{6{\kappa _{max}}}}{{{C_0}}}\rho + {2{\rho ^2}} .
\eeqnum
From here, we obtain a supplementary condition:
\beqnum\label{demo_9}
\frac{{2{\kappa _{max}}}}{{{C_0}}}\rho  < 1.
\eeqnum
This condition, along with \emph{Cond1} presented before, is equivalent to:\\
$$
(Cond2): 9\rho  < \frac{{{\kappa _{max}}}}{{{C_0}}} < \frac{1}{{2\rho }}.\\
$$
Therefore, $C_1$ has to be chosen, such that:
\beqnum\label{demo_10}
(Cond3):{C_1} > \frac{{\frac{{6{\kappa _{\max }}}}{{{C_0}}}\rho  + 2{\rho ^2}}}{{1 - \frac{{2\rho {\kappa _{\max }}}}{{{C_0}}}}}.\\
\eeqnum
\textbf{Case 2:} $\left| {{y_1}} \right| < \frac{1}{M}$.\\

As the saturation is no longer activated and $\left| \xi  \right| \leqslant 2\rho $, we obtain:
\beqnum
A \ge {C_1}My_1^2 - \left( {\frac{{\left( {3 + {C_1}} \right){\kappa _{\max }}}}{{{C_0}}} + 2\rho } \right)\left| {\xi {y_1}} \right| + \frac{1}{2}\left( {N - \frac{1}{{{C_0}}}} \right){\xi ^2}.\\
\eeqnum

$A$ is greater than quadratic form. To prove that it is positive definite, it is sufficient that
\beqnum\label{demo_11}
\begin{gathered}
  \left( {N - \frac{1}{{{C_0}}}} \right) > 0 \hfill \\
	\left| {\begin{array}{*{20}{c}}
   {{C_1}M} & { - \frac{{{\kappa _{max}}\left( {3 + {C_1}} \right)}}{{2{C_0}}} - \rho }  \\
   { - \frac{{{\kappa _{max}}\left( {3 + {C_1}} \right)}}{{2{C_0}}} - \rho } & {\frac{{N - \frac{1}{{{C_0}}}}}{2}}  \\
			\end{array}} \right| > 0.\hfill \\
\end{gathered}
\eeqnum
Equation (\ref{demo_11}) gives us:
\beqnum\label{demo_112}
{C_1}M\frac{{N - \frac{1}{{{C_0}}}}}{2} > {\left( {\frac{{{\kappa _{max}}\left( {3 + {C_1}} \right)}}{{2{C_0}}} + \rho } \right)^2}.\\
\eeqnum
Therefore, $M$ should be chosen such that:
\beqnum\label{demo_12}
(Cond4): M > \frac{{2{{\left( {\frac{{{\kappa _{max}}\left( {3 + {C_1}} \right)}}{{2{C_0}}} + \rho } \right)}^2}}}{{{C_1}\left( {N - \frac{1}{{{C_0}}}} \right)}}.\\
\eeqnum
\ede

\emph{Proof of Lemma \ref{lem_demo2}:} $B$ can be expressed in the following manner:\\
\beqnum\label{demo_13}
B = \left( {1 - \frac{{2{\rho ^2}}}{3}} \right)\left( {{y_2} - \frac{{\sigma \left( {{C_2}{y_2}} \right)}}{{{C_2}}}} \right)\sigma \left( {{C_2}{y_2}} \right) + \frac{\rho }{{{C_2}}}D\left( {\sigma \left( {{C_2}{y_2}} \right),\xi } \right),
\eeqnum
where
\beqnum\label{demo_14}
D\left( {z,\xi } \right): = \left( {1 - \frac{{2{\rho ^2}}}{3}} \right){z^2} - {C_2}N\left| {\xi z} \right| + \frac{{{C_2}}}{\rho }\left( {N - \frac{1}{{{C_0}}}} \right){\xi ^2}.
\eeqnum
It can be seen from equations (\ref{demo_13}) and (\ref{demo_14}), that if $D$ is positive definite, then $B$ is positive definite as well, i.e.:
\beqnum\label{demo_15}
\left| {\begin{array}{*{20}{c}}
  {1 - \frac{{2{\rho ^2}}}{3}}&{\frac{{ - {C_2 N}}}{2}} \\
  {\frac{{ - {C_2 N}}}{2}}&{\frac{{{C_2}}}{\rho }\left( {N - \frac{1}{{{C_0}}}} \right)}
\end{array}} \right| > 0.
\eeqnum
From here, we obtain a new condition on $\rho$:
\beqnum\label{demo_16}
(Cond5): \frac{{\left( {1 - \frac{{2{\rho ^2}}}{3}} \right)}}{\rho } > \frac{{{C_2}{N^2}}}{{4\left( {N - \frac{1}{{{C_0}}}} \right)}}.
\eeqnum
\ede

Therefore, to prove the theorem 1, it has to be shown that there exist constants $C_0$, $C_1$, $C_2$, $M$, $N$, $\rho$, such that conditions ($Cond2$) to ($cond5$) are met. In practice, $C_0$ and $C_2$ are given fixed positive values, and then $N$ is fixed such that $N > \frac{1}{{{C_0}}}$. Then, $\rho$ is chosen, small enough to satisfy conditions ($cond2$) and ($cond5$). Finally, $C_1$ and $M$ are chosen so that they satisfy respectively conditions ($cond3$) and ($cond4$).
\ede

The results presented above can be improved in the following directions

\begin{rem}\label{rem1}
Having a strict Lyapunov function allows us to extend the precedent results to cases in which external perturbations exist. More precisely, it can be shown that (\ref{system2}) is ISS (input-to-state) with respect to bounded external perturbations and an upper bound for allowed perturbations can be determined explicitly (as a function of the constants of the problem). In particular, it is interesting to suppose that the reference trajectory curvature $\kappa_r$, along with the vehicle velocity $V_x$ are susceptible to measurement noise. Hence the system can be stabilized in the proximity of the reference curve, depending explicitly on the magnitude of noise. In the following section, we will present simulation results, both with and without perturbations on $\kappa_r$.
\end{rem}




\begin{rem}\label{rem3}
It is possible not to bound the control $u_2$ as defined in (\ref{oup0}) but to simply use
$$
u_2=-C_0\left[\xi + \rho \sigma(C_2y_2)\right].
$$
The proof of the non-explosion of (\ref{edo-v}) is slightly modified but straigthforward.
\end{rem}



\section{Simulations}
In order to illustrate the performance of the presented controller, let us consider a unicycle type vehicle, with the following parameters:

$$
d=2\ m, \ V_x=15\ m.s^{-1}.
$$

The maximum curvature in the simulation is bounded by $\kappa _{max}=0.02 \ m^{-1}$. In order to highlight our claim that the performance of the controller is global and independent of initial condition, the value of $d\kappa_{max}$ has been kept much smaller than $1$ (in particular, $\xi(0)$ close to $\pi$).

The initial conditions imposed upon the error are
$$
e_p(0)=e_q(0)=10\ m,\  \xi(0)=9\pi/10.
$$

The control law $u_1,\ u_2$ are defined by
$$
u_1 = C_1 \sigma(M y_1),\\
u_2 = -\beta \sigma({\displaystyle\frac {C_0}{\beta}\left[\xi + \rho
\sigma(C_2 y_2)\right]}),
$$
where, the parameters have been determined according to \textit{Lemma 2.3} and \textit{Lemma 2.5}, specifically:
$$
C_0=0.4 \ ,\ C_1=0.7 \ ,\ C_2=1 \ ,\ M=1562 \ ,\ \beta=0.96 \ ,\ \rho=0.2 \ .
$$

\begin{figure}[h]
 \centering
      \includegraphics[width=9.5cm]{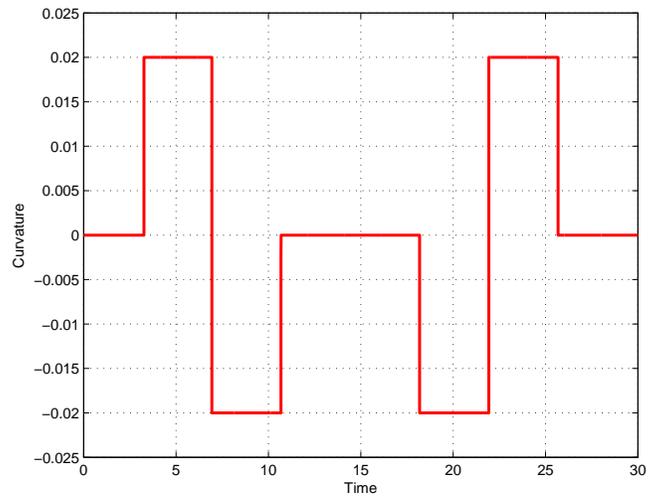}
 \caption{curvature $\kappa_r$ } \label{kappa}
\end{figure}

\begin{figure}[b]
 \centering
      \includegraphics[width=9.5cm]{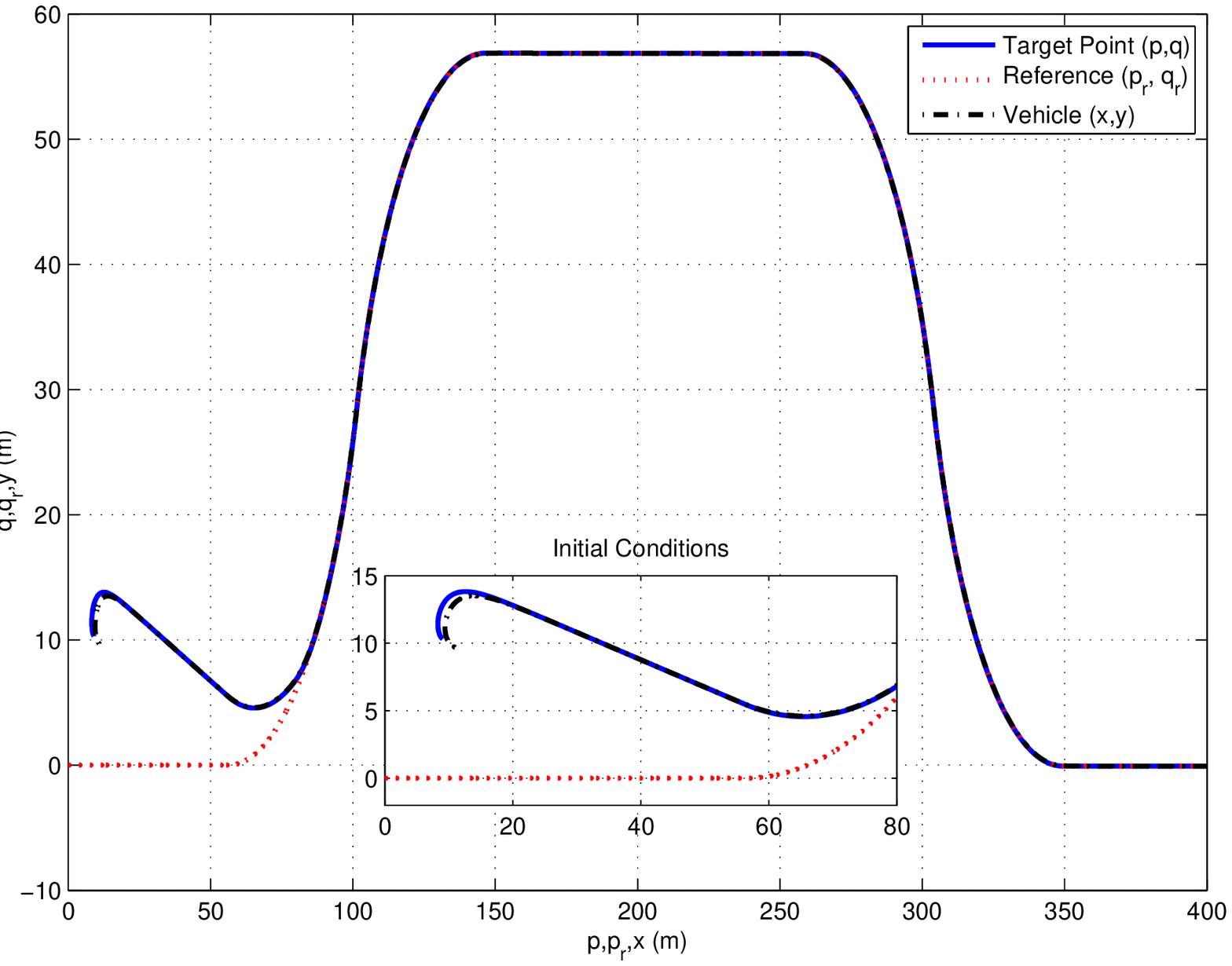}
  \caption{Reference trajectory, of the vehicle and its target point.} \label{tra1}
\end{figure}

\begin{figure}[t]\hspace{-1.0cm}
 \centering
      \includegraphics[width=9.5cm]{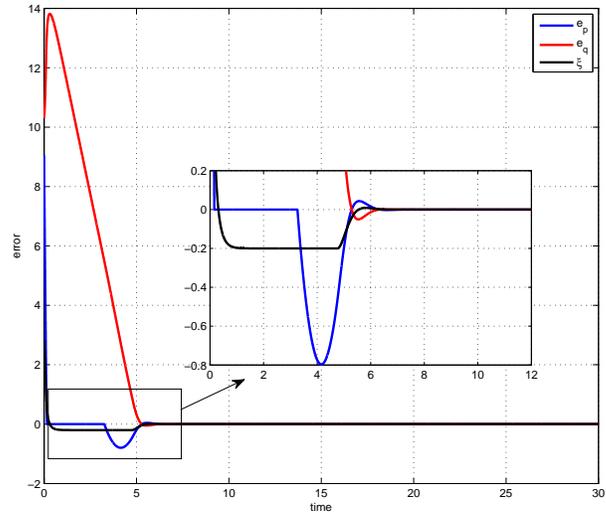}
  \caption{Errors $e_p$, $e_q$ et $\xi$} \label{erreurs}
\end{figure}

\begin{figure}[h]\hspace{-1.0cm}
 \centering
      \includegraphics[width=9.5cm]{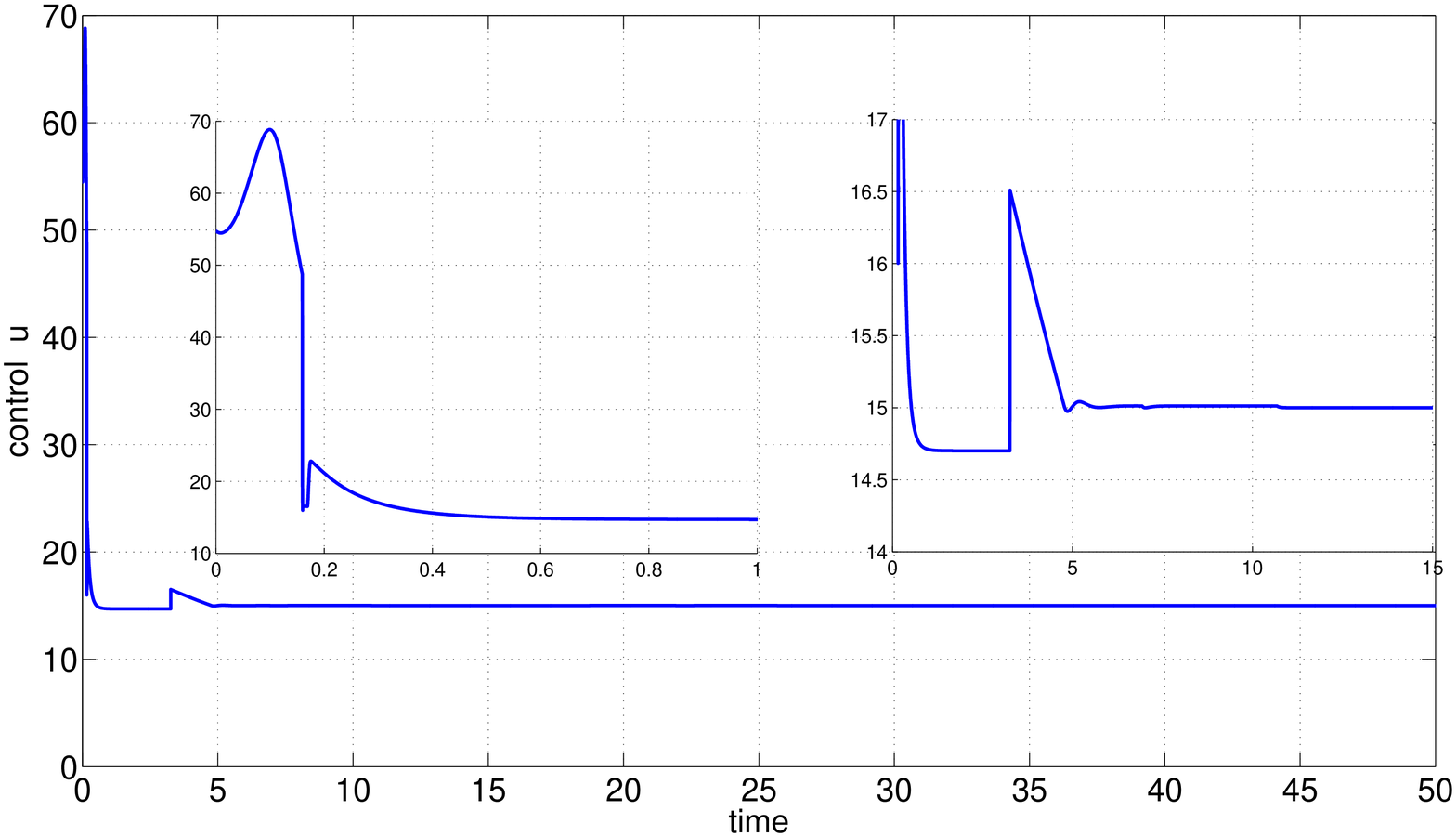}
  \caption{Control $u$} \label{com_u}
\end{figure}

\begin{figure}[b]\hspace{-1.0cm}
 \centering
      \includegraphics[width=9.5cm]{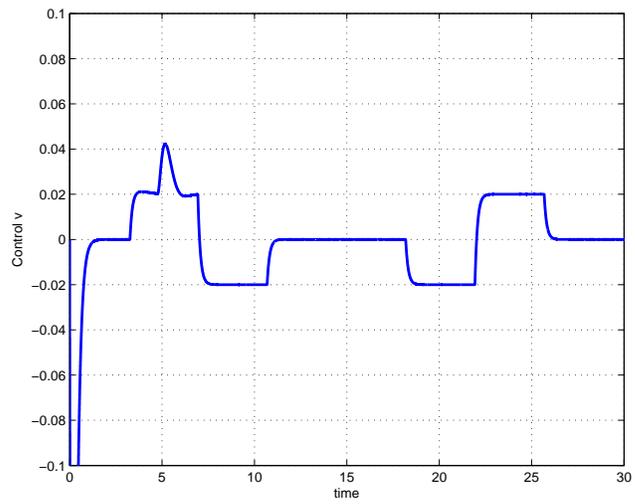}
  \caption{Control $v$} \label{com_v}
\end{figure}

The path to be followed $\gamma$ is defined by the geodesic curvature $\kappa_r$ (see Fig.~\ref{kappa}). It can be seen that the vehicle follows the target reference path as shown in Figure \ref{tra1}. The target point trajectory converges on the path in approximately $7$ sec. (see Fig.~\ref{erreurs}), and the vehicle successfully tracks the reference trajectory. The graphs of the control function are given in Fig.~\ref{com_u} and Fig.~\ref{com_v}.

\section{Conclusion}

In this article, we have addressed the problem of path following using a target point rigidly attached to a unicycle type vehicle. The control has been implemented using only the orientation of the vehicle. The main idea is to consider the parametrization of the followed path as an additional input for the system defined by the error dynamics. Control laws using saturation have been determined in order to achieve global stabilization without restrictions on initial conditions, under a (necessary) geometric characterization hypothesis, which relates the followed path with the target point position. This approach can also be extended to the cases where there are external perturbations or uncertainties in the model. This work can be extended towards addressing similar issues in more elaborate car models.

\bibliography{CIFA06}


\end{document}